\documentclass{siamltex}%
\usepackage{amsfonts}
\usepackage{amsmath}
\usepackage{amssymb}
\usepackage{graphicx}%
\newtheorem{remark}[theorem]{Remark}

\begin{document}

\title{On the Equivalence of Primal and Dual Substructuring
Preconditioners\thanks{This research was supported by the National Science
Foundation under grants \mbox{CNS-0325314}, \mbox{DMS-071387}, and
\mbox{CNS-0719641}.}}
\author{Bed\v{r}ich Soused\'ik\thanks{Department of Mathematics, Faculty of Civil
Engineering, Czech Technical University, Th\'akurova 7, 166 29 Prague 6, Czech
Republic and Department of Mathematical and Statistical Sciences, University of Colorado
Denver, P.O. Box 173364, Campus Box 170, Denver, CO 80217, USA
(\texttt{bedrich.sousedik@ucdenver.edu}). Supported in part by the program of
the Information society of the Academy of Sciences of the Czech Republic
\mbox{1ET400760509} and by the Grant Agency of the Czech Republic
\mbox{GA \v{C}R 106/08/0403}.}
\and Jan Mandel\thanks{Department of Mathematical and Statistical Sciences, University of Colorado
Denver, P.O. Box 173364, Campus Box 170, Denver, CO 80217, USA
(\texttt{jan.mandel@ucdenver.edu}).} }
\maketitle

\begin{abstract}
After a short historical review, we present four popular substructuring
methods: FETI-1, BDD, FETI-DP, BDDC, and derive the primal versions to the two
FETI methods, called P-FETI-1 and P-FETI-DP, as proposed by Fragakis and
Papadrakakis. The formulation of the BDDC method shows that it is the same as
P-FETI-DP and the same as a preconditioner introduced by Cros. We prove the
equality of eigenvalues of a particular case of the FETI-1 method and of the
BDD method by applying a recent abstract result by Fragakis.

\end{abstract}



\section{Introduction}

\label{sec:introduction}

Substructuring methods are among the most popular and widely used methods for
the solution of systems of linear algebraic equations obtained by finite
element discretization of second order elliptic problems.
This paper provides a review of recent results on the equivalence of several substructuring methods in
a common framework, complemented by some details not published previously.

We first give a brief review of the history of these methods (Section~\ref{sec:history}).
After introducing the basic concepts of substructuring (Section~\ref{sec:formulation}), we formulate the dual methods, FETI-1 and FETI-DP (Section~\ref{sec:feti-methods}), and derive their primal
versions, P-FETI-1 and P-FETI-DP, originally introduced in~\cite{Fragakis-2003-MHP}. However the derivation was  omitted in~\cite{Fragakis-2003-MHP}. Next,
we formulate the primal methods, BDD and BDDC (Section~\ref{sec:BDD-methods}).
Finally, we study connections between the methods in Section~\ref{sec:connections}.
We revisit our recent proof that the P-FETI-DP is in fact the same method
as the BDDC~\cite{Mandel-2007-BFM} and the preconditioner by Cros~\cite{Cros-2003-PSC}.
Next, we translate some of the abstract ideas from~\cite{Fragakis-2007-FDD,Fragakis-2003-MHP} into a framework usual in the domain decomposition literature.  We recall from~\cite{Fragakis-2003-MHP} that for a certain variant of FETI-1, the P-FETI-1 method is the same algorithm as BDD. Then we derive a
recent abstract result by Fragakis~\cite{Fragakis-2007-FDD} in this special case to show that the
eigenvalues of BDD and that particular version of FETI-1 are the same. It is
notable that this is the variant of FETI-1 devised to deal with difficult,
heterogeneous problems~\cite{Bhardwaj-2000-AFM}.

\section{Historical remarks}

\label{sec:history}

In this section, we provide a short overview of iterative substructuring, also
known as non-overlapping domain decomposition. Rather than attempting a
complete unbiased survey, our review centers on works connected to the
BDD\ and FETI\ theory by the second author and collaborators.

Consider a second order, selfadjoint, positive definite elliptic problem, such
as the Laplace equation or linearized elasticity, discretized by finite
elements with characteristic element size $h$. Given sufficient boundary
conditions, the global stiffness matrix is nonsingular, and its condition
number grows as $O\left(  h^{-2}\right)  $ for $h\to0$. However, if the domain
is divided into substructures consisting of disjoint unions of elements and
the interior degrees of freedom of each substructure are eliminated, the
resulting matrix on the boundary degrees of freedom has a condition number that grows
only as $O\left(  H^{-1}h^{-1}\right)  $, where $H>>h$ is the characteristic
size of the substructure. This fact has been known early on (Keyes and Gropp
\cite{Keyes-1987-CDD}); for a recent rigorous treatment, see Brenner
\cite{Brenner-1999-CNS}. The elimination of the interior degrees of freedom is
also called \emph{static condensation}, and the resulting reduced matrix is
called the \emph{Schur complement}. Because of the significant decrease of the
condition number, one can substantially accelerate iterative methods by
investing some work up front in the Choleski decomposition of the stiffness
matrix on the interior degrees of freedom and then just run back substitution
in each iteration. The finite element matrix is assembled separately in each
substructure. This process is called \emph{subassembly}. The elimination of
the interior degrees of freedom in each substructure can be done
independently, which is important for parallel computing:\ each substructure
can be assigned to an independent processor. The substructures are then
treated as large elements, with the Schur complements playing the role of the
local stiffness matrices of the substructures.
See~\cite{Keyes-1987-CDD,Smith-1996-DD} for more details.

The process just described is the background of \emph{primal iterative
substructuring methods}. Here, the condition that the values of degrees of
freedom common to several substructures coincide is enforced strongly, by
using a single variable to represent them. The improvement of the condition
number from $O\left(  h^{-2}\right)  $ to $O\left(  H^{-1}h^{-1}\right)  $,
straightforward implementation, and the potential for parallel computing
explain the early popularity of iterative substructuring methods
\cite{Keyes-1987-CDD}. However, further preconditioning is needed. Perhaps the
most basic preconditioner for the reduced problem is a diagonal one.
Preconditioning of a matrix by its diagonal helps to take out the dependence
on scaling and variation of coefficients and grid sizes. But the diagonal of
the Schur complement is expensive to obtain. It is usually better to avoid
computing the Schur complement explicitly and only use multiplication by the
reduced substructure matrices, which can be implemented by solving a
\emph{Dirichlet} problem on each substructure. Probing methods (Chan and
Mathew~\cite{Chan-1992-APT}) use such matrix-vector multiplication to estimate
the diagonal entries of the Schur complement.

In \emph{dual iterative substructuring methods}, also called FETI methods, the
condition that the values of degrees of freedom common to several
substructures coincide is enforced weakly, by Lagrange multipliers. The
original degrees of freedom are then eliminated, resulting in a~system for the
Lagrange multipliers, with the system operator consisting essentially of an
assembly of the inverses of the Schur complements. Multiplication by the
inverses of the Schur complements can be implemented by solving
a~\emph{Neumann} problem on each substructure. The assembly process is
modified to ensure that the Neumann problems are consistent, giving rise to
a~natural coarse problem. The system for the Lagrange multipliers is solved
again iteratively. This is the essence of the FETI method by Farhat and
Roux~\cite{Farhat-1991-MFE}, later called FETI-1. The condition number of the
FETI-1 method with diagonal preconditioning grows as $O\left(  h^{-1}\right)
$ and is bounded independently of the number of substructures (Farhat, Mandel,
and Roux~\cite{Farhat-1994-OCP}). For a small number of substructures, the
distribution of the eigenvalues of the iteration operator is clustered at
zero, resulting in superconvergence of conjugate gradients; however, for more
than a handful of substructures, the superconvergence is lost and the speed of
convergence is as predicted by the $O\left(  h^{-1}\right)  $ growth of the
condition number~\cite{Farhat-1994-OCP}.

For large problems and large number of substructures, \emph{asymptotically
optimal preconditioners} are needed. These preconditioners result typically in
condition number bounds of the form $O\left(  \log^{\alpha}\left(
1+H/h\right)  \right)  $ (the number $1$ is there only to avoid the value
$\log1=0$). In particular, the condition number is bounded independently of
the number of substructures and the bounds grow only slowly with the
substructure size. Such preconditioners require a \emph{coarse problem}, and
\emph{local preconditioning} that inverts approximately (but well enough) the
diagonal submatrices associated with segments of the interfaces between the
subtructures or the substructure matrices themselves. The role of the local
preconditioning is to slow down the growth of the condition number as
$h\rightarrow0$, while the role of the coarse problem is to provide global
exchange of information in order to bound the condition number independently
of the number of substructures. Many such asymptotically optimal primal
methods were designed in the 1980s and 1990s, e.g., Bramble, Pasciak, and
Schatz \cite{Bramble-1986-CPE,Bramble-1989-CPE}, Dryja~\cite{Dryja-1988-MDD},
Dryja, Smith, and Widlund~\cite{Dryja-1994-SAI}, Dryja and Widlund
\cite{Dryja-1995-SMN}, and Widlund~\cite{Widlund-1988-ISM}. However, those
algorithms require additional assumptions and information that may not be
readily available from finite element software, such as an explicit assumption
that the substructures form a coarse triangulation and that one can build
coarse linear functions from its vertices.

Practitioners desire methods that work algebraically with arbitrary
substructures, even if a theory may be available only in special cases
(first results on extending the theory to quite arbitrary substructures are
given in Dohrmann, Klawonn, and Widlund~\cite{Dohrmann-2008-DDL} and
Klawonn, Rheinbach, and Widlund~\cite{Klawonn-2008-AFA}).
They also prefer methods formulated in terms of the substructure matrices only, with minimal
additional information. In addition, the methods should be robust with respect
to various irregularities of the problem. Two such methods have emerged in
early 1990s:\ the Finite Element Tearing and Interconnecting (FETI)\ method by
Farhat and Roux~\cite{Farhat-1991-MFE}, and the Balancing Domain Decomposition
(BDD) by Mandel~\cite{Mandel-1993-BDD}. Essentially, the FETI method (with the
Dirichlet preconditioner) preconditions the assembly of the inverses of the Schur
complements by an assembly of the Schur complements, and the BDD method
preconditions assembly of Schur complements by an assembly of the inverses,
with a suitable coarse problem added. Of course, the assembly weights and
other details play an essential role.

The BDD method added a coarse problem to the local Neumann-Neumann
preconditioner by DeRoeck and Le Tallec \cite{DeRoeck-1991-ATL}, which
consisted of the assembly (with weights) of pseudoinverses of the local
matrices of the substructure. Assembling the inverses of the local
matrices is an idea similar to the Element-by-Element (EBE) method by Hughes
et al.~\cite{Hughes-1983-ESA}. The method was called Neumann-Neumann because
the preconditioner requires solution of Neumann problems on all substructures,
in contrast to an earlier Neumann-Dirichlet method, which, for a problem with
two substructures, required the solution of a Neumann problem on one and a
Dirichlet problem on the other~\cite{Widlund-1988-ISM}. The coarse problem in
BDD\ was constructed from the natural nullspace of the problem (constant for
the Laplace equation, rigid body motions for elasticity) and solving the
coarse problem guaranteed consistency of local problems in the preconditioner.
The coarse correction was then imposed variationally, just as the coarse
correction in multigrid methods. The $O\left(  \log^{2}\left(  1+H/h\right)
\right)  $ bound was then proved~\cite{Mandel-1993-BDD}.

In the FETI method, solving the local problems on the substructures to
eliminate the original degrees of freedom has likewise required working in the
complement of the nullspace of the substructure matrices, which gave a rise to
a natural coarse problem. Since the operator employs inverse of the Schur
complement (solving a Neumann problem) an optimal preconditioner employs
multiplication by the Schur complement (solving a Dirichlet problem), hence
the preconditioner was called the Dirichlet preconditioner. The $O\left(
\log^{3}\left(  1+H/h\right)  \right)  $ bound was proved by Mandel and
Tezaur~\cite{Mandel-1996-CSM}, and $O\left(  \log^{2}\left(  1+H/h\right)
\right)  $ for a certain variant of the method by
Tezaur~\cite{Tezaur-1998-ALM}. See also Klawonn and Widlund~\cite{Klawonn-2001-FNN}
for further discussion.

Because the interface to the BDD\ and FETI method required only the
multiplication by the substructure Schur complements, solving systems with the
substructure Schur complements, and information about the substructure
nullspace, the methods got quite popular and widely used. In Cowsar, Mandel,
and Wheeler~\cite{Cowsar-1995-BDD}, the multiplications were implemented as
solution of mixed problems on substructures. However, neither the BDD nor the
FETI method worked well for 4th order problems (plate bending). The reason was
essentially that both methods involve \textquotedblleft
tearing\textquotedblright\ a~vector of degrees of freedom reduced to the
interface, and, for 4th order problems, the \textquotedblleft
torn\textquotedblright\ function has energy that grows as negative power of
$h$, unlike for 2nd order problems, where the energy grows only as a positive
power of $\log1/h$. The solution was to prevent the \textquotedblleft
tearing\textquotedblright\ by fixing the function at the substructure corners;
then only its derivative along the interface gets \textquotedblleft
torn\textquotedblright, which has energy again only of the order $\log1/h$.
Preventing such \textquotedblleft tearing\textquotedblright\ can be generally
accomplished by increasing the coarse space, since the method runs in the
complement to the coarse space. For the BDD method, this was relatively
straightforward, because the algebra of the BDD method allows arbitrary
enlargement of the coarse space. The coarse space that does the trick contains
additional functions with spikes at corners, defined by fixing the value at
the corner and minimizing the energy. With this improvement, $O\left(
\log^{2}\left(  1+H/h\right)  \right)  $ condition number bound was proved and
fast convergence was recovered for 4th order problems (Le Tallec, Mandel, and
Vidrascu~\cite{LeTallec-1994-BDD,LeTallec-1998-NND}). In the FETI\ method,
unfortunately, the algebra requires that the coarse space is made of exactly
the nullspace of the substructure matrices, so a simple enlargement of the
coarse space is not possible. Therefore, a version of FETI, called FETI-2, was
developed by Mandel, Tezaur, and Farhat~\cite{Mandel-1999-SSM}, with a second
correction by coarse functions concentrated at corners, wrapped around the
original FETI method variationally much like BDD, and the $O\left(  \log
^{3}\left(  1+H/h\right)  \right)  $ bound was proved again. However, the BDD
and FETI methods with the modifications for 4th order problems were rather
unwieldy (especially FETI-2), and, consequently, not as widely used.

The breakthrough came with the Finite Element Tearing and Interconnecting - Dual, Primal (FETI-DP) method by Farhat et
al.~\cite{Farhat-2001-FDP}, which enforced the continuity of the degrees of
freedom on a substructure corner as in the primal method by representing them
by one common variable, while the remaining continuity conditions between the
substructures are enforced by Lagrange multipliers. The primal variables are
again eliminated and the iterations run on the Lagrange multipliers. The
elimination process can be organized as solution of sparse system and it gives
rise to a natural coarse problem, associated with substructure corners. In 2D,
the FETI-DP method was proved to have condition number bounded as $O\left(
\log^{2}\left(  1+H/h\right)  \right)  $ both for 2nd order and 4th order
problems by Mandel and Tezaur~\cite{Mandel-2001-CDP}. However, the method does
not converge as well in 3D and averages over edges or faces of substructures
need to be added as coarse variables for fast convergence (Klawonn, Widlund,
and Dryja~\cite{Klawonn-2002-DPF}, Farhat, Lesoinne, and Pierson
\cite{Farhat-2000-SDP}), and the $O\left(  \log^{2}\left(  1+H/h\right)
\right)  $ bound can then be proved~\cite{Klawonn-2002-DPF}.

The Balancing Domain Decomposition by Constraints (BDDC) was developed by
Dohrmann~\cite{Dohrmann-2003-PSC} as a~primal alternative the
FETI-DP\ method. The BDDC\ method uses imposes the equality of coarse degrees
of freedom on corners and of averages by constraints. In the case of only
corner constraints, the coarse basis functions are the same as in the BDD
method for 4th order problems from~\cite{LeTallec-1994-BDD,LeTallec-1998-NND}.
The bound $O\left(  \log^{2}\left(  1+H/h\right)  \right)  $ for BDDC was
first proved by Mandel and Dohrmann~\cite{Mandel-2003-CBD}. The BDDC\ and the
FETI-DP\ are currently the most advanced versions of the BDD and FETI families
of methods.

The convergence properties of the BDDC and FETI-DP methods were quite similar,
yet it came as a surprise when Mandel, Dohrmann, and Tezaur
\cite{Mandel-2005-ATP} proved that the spectra of their preconditioned
operators are in fact identical, once all the components are same. This result
came at the end of a long chain of ties discovered between BDD and FETI\ type
method. Algebraic relations between FETI and BDD methods were pointed out by
Rixen et al.~\cite{Rixen-1999-TCF}, Klawonn and
Widlund~\cite{Klawonn-2001-FNN}, and Fragakis and
Papadrakakis~\cite{Fragakis-2003-MHP}. An important common bound on the
condition number of both the FETI and the BDD method in terms of a single
inequality was given Klawonn and Widlund~\cite{Klawonn-2001-FNN}. Fragakis and
Papadrakakis~\cite{Fragakis-2003-MHP}, who derived certain primal versions of
FETI and FETI-DP preconditioners (called P-FETI-1 and P-FETI-DP), have also
observed that the eigenvalues of BDD and a certain version of FETI are
identical along with the proof that the primal version of this particular FETI
algorithm gives a method same as BDD. The proof of equality of eigenvalues of
BDD and FETI was given just recently in more abstract framework by
Fragakis~\cite{Fragakis-2007-FDD}. Mandel, Dohrmann, and
Tezaur~\cite{Mandel-2005-ATP} have proved that the eigenvalues of BDDC and
FETI-DP are identical and they have obtained a simplified and fully algebraic
version (i.e., with no undetermined constants) of a common condition number
estimate for BDDC and FETI-DP, similar to the estimate by Klawonn and
Widlund~\cite{Klawonn-2001-FNN} for BDD and FETI. Simpler proofs of the
equality of eigenvalues of BDDC and FETI-DP were obtained by Li and
Widlund~\cite{Li-2006-FBB}, and by Brenner and Sung~\cite{Brenner-2007-BFW},
who also gave an example when BDDC\ has an eigenvalue equal to one but
FETI-DP\ does not. A~primal variant of P-FETI-DP was proposed by
Cros~\cite{Cros-2003-PSC}, giving a conjecture that P-FETI-DP\ and BDDC\ is in
fact the same method, which was first shown on a somehow more abstract level in our recent
work~\cite{Mandel-2007-BFM}.

It is interesting to note that the choice of assembly weights in the BDD
preconditioner was known at the very start from the work of DeRoeck and Le
Tallec \cite{DeRoeck-1991-ATL} and before, while the choice of weights for
FETI type method is much more complicated. A correct choice of weights is
essential for the robustness of the methods with respect to scaling the matrix
in each substructure by an arbitrary positive number (the \textquotedblleft
independence of the bounds on jumps in coefficients\textquotedblright). For
the BDD method, such convergence bounds were proved~by Mandel and Brezina
\cite{Mandel-1996-BDD}, using a similar argument as in
Sarkis~\cite{Sarkis-1993-TSM} for Schwarz methods; see also Dryja, Sarkis, and
Widlund~\cite{Dryja-1996-MSM}. For the FETI methods, a proper choice of
weights was discovered only much later - see Rixen and Farhat~\cite{Rixen-1999-SEE}, Farhat, Lesoinne and
Pierson~\cite{Farhat-2000-SDP} for a special cases, Klawonn and
Widlund~\cite{Klawonn-2001-FNN} for a more general case and convergence
bounds, and a detailed discussion in Mandel, Dohrmann, and
Tezaur~\cite{Mandel-2005-ATP}.

\section{Substructuring Components for a Model Problem}

\label{sec:formulation}We first show how the spaces and operators we will work
with arise in the standard substructuring theory for a model problem obtained
by a discretization of the second order elliptic problem. Consider a bounded
domain $\Omega\subset\mathbb{R}^{d}$ decomposed into nonoverlapping subdomains
(alternatively called substructures) denoted $\Omega_{i}$, $i=1,...,N$, which
form a conforming triangulation of the domain$~\Omega$. Each substructure is a
union of a uniformly bounded number of Lagrangean $P1$ or $Q1$ finite
elements, such that the nodes of the finite elements between substructures
coincide. The boundary of$~\Omega_{i}$ is denoted by $\partial\Omega_{i}$. The
nodes contained in the intersection of at least two substructures are called
boundary nodes. The union of all boundary nodes of all substructures is called
the interface $\Gamma$ and\ $\Gamma_{i}$ is the interface of
substructure$~\Omega_{i}$. The space of vectors of local degrees of freedom
on$~\Gamma_{i}$ is denoted by $W_{i}$ and $W=W_{1}\times\cdots\times W_{N}$.
Let $S_{i}:W_{i}\rightarrow W_{i}$ be the Schur complement operator obtained
by eliminating all interior degrees of freedom of $\Omega_{i}$, i.e., those
that do not belong to interface $\Gamma_{i}$. We assume that the matrices
$S_{i}$ are symmetric positive semidefinite and consider global vectors and
matrices in the block form
\begin{equation}
w=\left[
\begin{array}
[c]{c}%
w_{1}\\
\vdots\\
w_{N}%
\end{array}
\right]  ,\quad w\in W,\quad S=
\begin{bmatrix}
S_{1} &  & \\
& \ddots & \\
&  & S_{N}%
\end{bmatrix}
. \label{eq:block-operators}%
\end{equation}
The problem we wish to solve is the constrained minimization of energy,
\begin{equation}
\frac{1}{2}a\left(  u,u\right)  -\left\langle r,u\right\rangle \rightarrow
\min\text{ subject to }u\in\widehat{W}, \label{eq:constrained-minim}%
\end{equation}
where $\widehat{W}\subset W$ is the space of all vectors of degrees of freedom
on the substructures that coincide on the interfaces, and the bilinear form
\[
a\left(  u,v\right)  =\left\langle Su,v\right\rangle ,\quad\forall u,v\in W,
\]
is assumed to be positive definite on $\widehat{W}$. In the variational form,
problem (\ref{eq:constrained-minim}) can be written as
\begin{equation}
u\in\widehat{W}:a(u,v)=\left\langle r,v\right\rangle ,\quad\forall
v\in\widehat{W}. \label{eq:problem-reduced}%
\end{equation}
The global Schur complement $\widehat{S}:\widehat{W}\mapsto\widehat{W}%
^{\prime}$ associated with $a$ is defined by
\begin{equation}
a(u,v)=\left\langle \widehat{S}u,v\right\rangle ,\quad\forall u,v\in
\widehat{W}. \label{eq:S-hat}%
\end{equation}
Defining $R$ as the natural embedding of the space $\widehat{W}$ into the
space $W$, i.e.,
\begin{equation}
R:\widehat{W}\rightarrow W,\quad R:u\in\widehat{W}\longmapsto u\in W,
\label{eq:def-r}%
\end{equation}
we can write (\ref{eq:problem-reduced}) equivalently as the system of linear
algebraic equations
\begin{equation}
\widehat{S}u=r,\quad\text{ where }\widehat{S}=R^{T}SR. \label{eq:system}%
\end{equation}

The BDDC\ and FETI-DP,\ as the two-level preconditioners, are characterized by
the selection of certain \emph{coarse degrees of freedom}, such as values at
the corners and averages over edges or faces of substructures (for their
general definition see, e.g.,~\cite{Klawonn-2006-DPF}). So, we define
$\widetilde{W}\subset W$ as the subspace of all functions such that the values
of any coarse degrees of freedom have a common value for all relevant
substructures and vanish on $\partial\Omega$, and such that%
\[
\widehat{W}\subset\widetilde{W}\subset W.
\]
The space $\widetilde{W}$ has to be selected in the design of the
preconditioner so that the bilinear form $a\left(  \cdot,\cdot\right)  $ is
positive definite on $\widetilde{W}$. The operator $\widetilde{S}%
:\widetilde{W}\rightarrow\widetilde{W}^{\prime}$ associated with $a$ is
defined by
\[
a(u,v)=\left\langle \widetilde{S}u,v\right\rangle ,\quad\forall u,v\in
\widetilde{W}.
\]

\begin{remark}
The idea to restrict the bilinear form $a(\cdot,\cdot)$ from the space $W$
into the subspace $\widetilde{W}$ is closely related to the concept of
subassembly as employed in~\cite{Li-2006-FBB}.
\end{remark}

In formulation of dual methods from the FETI\ family, we introduce the matrix%
\[
B=\left[  B_{1},\ldots,B_{N}\right]  ,
\]
which enforces the continuity across\ substructure interfaces and it is
defined as follows: each row $B$ corresponds to a degree of freedom common to
a pair of substructures $i$ and $j$. The entries of the row are zero except
for one $+1$ in the block $i$ and one $-1$ in the block $j$, so that the
condition
\[
Bu=0\Longleftrightarrow u\in\widehat{W},
\]
and using (\ref{eq:def-r}), clearly
\begin{equation} BR=0. \label{eq:BR} \end{equation}
An important ingredient of substructuring methods is the averaging operator
$E:W\rightarrow\widehat{W}$\ defined as%
\begin{equation}
E=R^{T}D_{P}, \label{eq:def-averaging}%
\end{equation}
where $D_{P}:W\rightarrow W$ is a given weight matrix such that the
decomposition of unity property holds,%
\begin{equation}
ER=I. \label{eq:decomp-unity}%
\end{equation}
In terms of substructuring, $E$ is an averaging operator that maps the
substructure local degrees of freedom to global degrees of freedom.

The last ingredient is the matrix $B_{D}$ constructed from $B$ as\
\[
B_{D}=\left[  D_{D1}B_{1},\ldots,D_{DN}B_{N}\right]  ,
\]
where the matrices $D_{Di}$ are determined from $D_{P}$, see
\cite{Klawonn-2002-DPF,Mandel-2005-ATP} for details.

Finally, we shall assume, cf., e.g., \cite[eq. (10)]
{Mandel-2005-ATP}, that
\begin{equation}
B_{D}^{T}B+RE =I,\label{eq:jump+average}
\end{equation}
which easily implies $EB_D^TB=E(I-RE)=E - ERE=0$, and so
\begin{equation} B^{T}B_{D}E^{T} =0. \label{eq:B_DE^T}%
\end{equation}

\section{P-FETI Family of Methods}

We review the FETI-1 and FETI-DP preconditioners followed in each case by a
formulation of their primal versions denoted as P-FETI-1 and P-FETI-DP, respectively.

\label{sec:feti-methods}

\subsection{P-FETI-1}

In the case of the FETI-1 method, the problem~(\ref{eq:constrained-minim}) is
formulated as minimization of total subdomain energy subject to the continuity
condition%
\begin{equation}
\frac{1}{2}a\left(  w,w\right)  -\left\langle f,w\right\rangle \rightarrow
\min\quad\ \text{subject to }w\in W,\quad Bw=0,
\label{eq:constraint-minim-feti}%
\end{equation}
which is equivalent to a saddle point system: find $\left(  w,\lambda\right)
\in W\times\Lambda$ such that
\begin{equation}%
\begin{array}
[c]{ccccc}%
Sw & + & B^{T}\lambda & = & f,\\
Bw &  &  & = & 0.
\end{array}
\label{eq:feti-system}%
\end{equation}
First, note that $S$ is invertible on $\operatorname*{null}B$ and $\lambda$ is
unique up to a component in $\operatorname*{null}B^{T}$, so $\Lambda$ is
selected to be $\operatorname*{range}B$. Let $Z$ be matrix with linearly
independent columns, such that
\begin{equation}
\operatorname*{range}Z=\operatorname*{null}S. \label{eq:rangeZ=nullS}%
\end{equation}
Since $S$ is semi-definite, it must hold for the first equation to be solvable
that
\[
f-B^{T}\lambda\in\operatorname*{range}S=\left(  \operatorname*{null}%
S\right)  ^{\perp}=\left(  \operatorname*{range}Z\right)  ^{\perp
}=\operatorname*{null}Z^{T},
\]
so, equivalently, we require that
\begin{equation}
Z^{T}(f-B^{T}\lambda)=0. \label{eq:nullZ^T}%
\end{equation}
Eliminating $w$ from the first equation of (\ref{eq:feti-system}) as
\begin{equation}
w=S^{+}(f-B^{T}\lambda)+Za, \label{eq:feti-w}%
\end{equation}
substituting in the second equation of (\ref{eq:feti-system}) and rewriting
(\ref{eq:nullZ^T}), we get
\[%
\begin{array}
[c]{ccccc}%
BS^{+}B^{T}\lambda & - & BZa & = & BS^{+}f,\\
-Z^{T}B^{T}\lambda &  &  & = & -Z^{T}f.
\end{array}
\]
Denoting $G=BZ$ and $F=BS^{+}B^{T}$ this system becomes
\begin{equation}%
\begin{array}
[c]{ccccc}%
F\lambda & - & Ga & = & BS^{+}f,\\
-G^{T}\lambda &  &  & = & -Z^{T}f.
\end{array}
\label{eq:feti-system-1}%
\end{equation}
Multiplying the first equation by $\left(  G^{T}QG\right)  ^{-1}G^{T}Q$, where
$Q$ is some symmetric and positive definite scaling\ matrix, we can compute
$a$ as
\begin{equation}
a=\left(  G^{T}QG\right)  ^{-1}G^{T}Q(F\lambda - BS^{+}f). \label{eq:feti-a}%
\end{equation}
The first equation in (\ref{eq:feti-system-1}) thus becomes%
\begin{equation}
F\lambda-G\left(  G^{T}QG\right)  ^{-1}G^{T}Q( F\lambda - BS^{+}f)=BS^{+}f.
\label{eq:feti-system-P^T}%
\end{equation}
Introducing
\begin{equation}
P=I-QG(G^{T}QG)^{-1}G^{T}, \label{eq:P}%
\end{equation}
as the $Q$-orthogonal projection onto $\operatorname*{null}G^{T}$, we get that
(\ref{eq:feti-system-P^T}) corresponds to the first equation in
(\ref{eq:feti-system-1}) multiplied by $P^{T}$. So, the system
(\ref{eq:feti-system-1}) can be written in the decoupled form as%
\[%
\begin{array}
[c]{ccc}%
P^{T}F\lambda & = & P^{T}BS^{+}f,\\
G^{T}\lambda & = & Z^{T}f.
\end{array}
\]
The initial value of $\lambda$ is chosen to satisfy the second equation in
(\ref{eq:feti-system-1}), so
\begin{equation}
\lambda_{0}=QG(G^{T}QG)^{-1}Z^{T}f. \label{eq:feti-lambda_0}%
\end{equation}
Substituting $\lambda_{0}$ into (\ref{eq:feti-a}) gives initial value of $a$
as%
\begin{equation}
a_{0}=\left(  G^{T}QG\right)  ^{-1}G^{T}Q(F\lambda_{0} - BS^{+}f).
\label{eq:feti-a_0}%
\end{equation}
Since we are looking for $\lambda\in\operatorname*{null}G^{T}$, the
FETI-1\ method is a preconditioned conjugate gradient method applied to the
system%
\begin{equation}
P^{T}FP\lambda=P^{T}BS^{+}f. \label{eq:feti}%
\end{equation}
In the primal version of the FETI-1 preconditioner, the assembled and averaged
solution $u$ is obtained from (\ref{eq:feti-w}), using equations
(\ref{eq:feti-a_0}) and (\ref{eq:feti-lambda_0}), as
\begin{align*}
u  &  =Ew\\
&  =E\left[  S^{+}(f-B^{T}\lambda_{0})+Za_{0}\right] \\
&  =E\left[  S^{+}(f-B^{T}\lambda_{0})+Z\left(  G^{T}QG\right)  ^{-1}%
G^{T}Q(F\lambda_{0}-BS^{+}f)\right] \\
&  =E\left[  S^{+}(f-B^{T}\lambda_{0})+Z\left(  G^{T}QG\right)  ^{-1}%
G^{T}Q(BS^{+}B^{T}\lambda_{0}-BS^{+}f)\right] \\
&  =E\left[  I-Z\left(  G^{T}QG\right)  ^{-1}G^{T}QB\right]  S^{+}%
(f-B^{T}\lambda_{0})\\
&  =E\left[  \left(  I-Z\left(  G^{T}QG\right)  ^{-1}G^{T}QB\right)
S^{+}\left(  I-B^{T}QG\left(  G^{T}QG\right)  ^{-1}Z^{T}\right)  \right]
E^{T}r\\
&  =EH^{T}S^{+}HE^{T}r\\
&  =M_{P-FETI}r,
\end{align*}
where we have denoted by
\begin{equation}
H=I-B^{T}QG\left(  G^{T}QG\right)  ^{-1}Z^{T}, \label{eq:H}%
\end{equation}
and so
\begin{equation}
M_{P-FETI}=EH^{T}S^{+}HE^{T}, \label{eq:P-FETI-1}%
\end{equation}
is the associated primal preconditioner P-FETI-1, same as~\cite[eq.
(79)]{Fragakis-2003-MHP}.

\subsection{P-FETI-DP}

In the case of the FETI-DP, the problem~(\ref{eq:constrained-minim}) is
formulated as minimization of total subdomain energy subject to the continuity
condition%
\begin{equation}
\frac{1}{2}a\left(  w,w\right)  -\left\langle f,w\right\rangle \rightarrow
\min\quad\ \text{subject to }w\in\widetilde{W},\quad Bw=0.
\label{eq:constraint-minim-feti-dp}%
\end{equation}
Compared to the formulation of FETI-1 in (\ref{eq:constraint-minim-feti}), we
have now used the subspace $\widetilde{W}\subset W$ such that the operator
$\widetilde{S}$ associated with $a\left(  \cdot,\cdot\right)  $ on the space
$\widetilde{W}$ is positive-definite. In this case,
(\ref{eq:constraint-minim-feti-dp}) is equivalent to setting up a saddle point
system: find $\left(  w,\lambda\right)  \in\widetilde{W}\times\Lambda$ such
that
\begin{equation}%
\begin{array}
[c]{ccccc}%
\widetilde{S}w & + & B^{T}\lambda & = & f,\\
Bw &  &  & = & 0.
\end{array}
\label{eq:FETI-DP-system}%
\end{equation}
Since $\widetilde{S}$ is invertible on $\widetilde{W}$, solving for $w$ from
the first and substituting into the second equation of
(\ref{eq:FETI-DP-system}), we get%
\begin{equation}
B\widetilde{S}^{-1}B^{T}\lambda=B\widetilde{S}^{-1}f,
\label{eq:FETI-DP-system2}%
\end{equation}
which is the dual system to be solved by preconditioned conjugate gradients,
with the Dirichlet preconditioner defined by
\begin{equation}
M_{FETI-DP}=B_{D}\widetilde{S}B_{D}^{T}. \label{eq:Dirichlet-preconditioner}%
\end{equation}

Next, we will derive the P-FETI-DP\ preconditioner using the original paper by
Farhat et. al.~\cite{Farhat-2001-FDP} in order to verify the P-FETI-DP
algorithm given in~\cite[eq. (90)]{Fragakis-2003-MHP} for the corner
constraints. We split the global vector of degrees of freedom $u$\ into the
vector of global coarse degrees of freedom denoted by $u_{c}$ and the vector
of remaining degrees of freedom denoted by $u_{r}$. We note that
we could perform a change of basis, cf., e.g.,~\cite{Klawonn-2001-DDP,Klawonn-2006-DPF,Li-2006-FBB}
to make all primal constraint (such as averages over edges or faces)
explicit, i.e., each coarse degrees of freedom would correspond to an explicit
degree of freedom in the vector $u_{c}$. Thus, we decompose the space
$\widetilde{W}$ as, cf.~\cite[Remark 5]{Mandel-2005-ATP},
\begin{equation}
\widetilde{W}=\widetilde{W}_{c}\oplus\widetilde{W}_{r},
\label{eq:tildeW-decomp}%
\end{equation}
where the space $\widetilde{W}_{c}$ consists of functions that are continuous
across interfaces, have a nonzero value at one coarse degree of freedom at a
time and zero at other coarse degrees of freedom, and the space $\widetilde
{W}_{r}$ consists of functions with coarse degrees of freedom equal to zero.
The solution splits into the solution of the global coarse problem in the
space $\widetilde{W}_{c}$ and the solution of independent subdomain problems
on the space $\widetilde{W}_{r}$.

Let $R_{c}^{(i)}$ be a map of global coarse variables to its subdomain
component, i.e.,
\[
R_{c}^{(i)}u_{c}=u_{c}^{(i)},\qquad R_{c}=\left(
\begin{array}
[c]{c}%
R_{c}^{(1)}\\
\vdots\\
R_{c}^{(N)}%
\end{array}
\right)  ,
\]
let $B_{r}$\ be an operator enforcing the interface continuity of $u_{r}$\ by
\[
B_{r}u_{r}=0,\qquad B_{r}=\left(
\begin{array}
[c]{ccc}%
B_{r}^{(1)} & \ldots & B_{r}^{(N)}%
\end{array}
\right)  ,
\]
and let the mappings $E_{r}^{T}$\ and $E_{c}^{T}$\ distribute the primal
residual $r$ to the subdomain forces and to the global coarse problem
right-hand side, respectively.

The equations of equilibrium can now be written, cf.~\cite[eq. (9)-(10)]%
{Farhat-2001-FDP}, as
\[%
\begin{array}
[c]{ccccccc}%
S_{rr}^{(i)}w_{r}^{(i)} & + & S_{rc}^{(i)}R_{c}^{(i)}w_{c} & + & B_{r}^{(i)T}\lambda & = & f_{r}^{(i)}, \\
\sum\limits_{i=1}^{N}R_{c}^{(i)T}S_{rc}^{(i)T}w_{r}^{(i)} & + &
\sum\limits_{i=1}^{N}R_{c}^{(i)T}S_{cc}^{(i)}R_{c}^{(i)}w_{c} &  &  & = &
f_{c}, \\
\sum\limits_{i=1}^{N}B_{r}^{(i)}w_{r}^{(i)} &  &  &  &  & = & 0,
\end{array}
\]
where the first equation corresponds to independent subdomain problems, second
corresponds to the global coarse problem and the third enforces the continuity
of local problems.\ This system can be re-written as%
\begin{equation}
\left(
\begin{array}
[c]{ccc}%
S_{rr} & S_{rc}R_{c} & B_{r}^{T}\\
\left(  S_{rc}R_{c}\right)  ^{T} & \widetilde{S}_{cc} & 0\\
B_{r} & 0 & 0
\end{array}
\right)  \left(
\begin{array}
[c]{c}%
u_{r}\\
u_{c}\\
\lambda
\end{array}
\right)  =\left(
\begin{array}
[c]{c}%
f_{r}\\
f_{c}\\
0
\end{array}
\right)  , \label{eq:FETI-DP-system-3}%
\end{equation}
where $f_{r}=E_{r}^{T}r$, $f_{c}=E_{c}^{T}r$,\ and the blocks are defined\ as
\[
\widetilde{S}_{cc}=\sum_{i=1}^{N}R_{c}^{(i)T}S_{cc}^{(i)}R_{c}^{(i)},\quad
S_{rr}=\left(
\begin{array}
[c]{ccc}%
S_{rr}^{(1)} &  & \\
& \ddots & \\
&  & S_{rr}^{(N)}%
\end{array}
\right)  ,\quad S_{rc}R_{c}=\left(
\begin{array}
[c]{c}%
S_{rc}^{(1)}R_{c}^{(1)}\\
\vdots\\
S_{rc}^{(N)}R_{c}^{(N)}%
\end{array}
\right)  .
\]

\begin{remark}
Note that the system (\ref{eq:FETI-DP-system-3}) is just the expanded system
(\ref{eq:FETI-DP-system}).
\end{remark}

Expressing $u_{r}$ from the first equation in (\ref{eq:FETI-DP-system-3}), we
get
\[
u_{r}=S_{rr}^{-1}\left(  f_{r}-S_{rc}R_{c}u_{c}-B_{r}^{T}\lambda\right)  .
\]
Substituting for $u_{r}$ into the second equation in
(\ref{eq:FETI-DP-system-3}) gives
\[
\widetilde{S}_{cc}^{\ast}u_{c}-\left(  S_{rc}R_{c}\right)  ^{T}S_{rr}%
^{-1}B_{r}^{T}\lambda=f_{c}-\left(  S_{rc}R_{c}\right)  ^{T}S_{rr}^{-1}f_{r},
\]
where $\widetilde{S}_{cc}^{\ast}=\widetilde{S}_{cc}-R_{c}^{T}S_{rc}^{T}%
S_{rr}^{-1}S_{rc}R_{c}$. Inverting $\widetilde{S}_{cc}^{\ast}$, we get that%
\[
u_{c}=\widetilde{S}_{cc}^{\ast^{-1}}\left[  f_{c}-\left(  S_{rc}R_{c}\right)
^{T}S_{rr}^{-1}f_{r}+\left(  S_{rc}R_{c}\right)  ^{T}S_{rr}^{-1}B_{r}%
^{T}\lambda\right]  .
\]
After initialization with $\lambda=0$, which
\cite{Fragakis-2003-MHP,Fragakis-2007-FDD} does not say, but it can be used,
cf., e.g.,~\cite[Section 6.4]{Toselli-2005-DDM}, the assembled and averaged
solution is
\begin{align*}
u  &  =E_{r}u_{r}+E_{c}u_{c}\\
&  =E_{r}S_{rr}^{-1}\left\{  f_{r}-S_{rc}R_{c}\widetilde{S}_{cc}^{\ast^{-1}%
}\left(  f_{c}-\left(  S_{rc}R_{c}\right)  ^{T}S_{rr}^{-1}f_{r}\right)
\right\}  +\\
&  \quad+E_{c}\widetilde{S}_{cc}^{\ast^{-1}}\left(  f_{c}-\left(  S_{rc}%
R_{c}\right)  ^{T}S_{rr}^{-1}f_{r}\right) \\
&  =E_{r}S_{rr}^{-1}f_{r}-E_{r}S_{rr}^{-1}S_{rc}R_{c}\widetilde{S}_{cc}%
^{\ast^{-1}}f_{c}+\\
&  \quad+E_{r}S_{rr}^{-1}S_{rc}R_{c}\widetilde{S}_{cc}^{\ast^{-1}}\left(
S_{rc}R_{c}\right)  ^{T}S_{rr}^{-1}f_{r}+\\
&  \quad+E_{c}\widetilde{S}_{cc}^{\ast^{-1}}f_{c}-E_{c}\widetilde{S}%
_{cc}^{\ast^{-1}}\left(  S_{rc}R_{c}\right)  ^{T}S_{rr}^{-1}f_{r}\\
&  =E_{r}S_{rr}^{-1}f_{r}+\\
&  \quad+\left(  E_{c}-E_{r}S_{rr}^{-1}S_{rc}R_{c}\right)  \widetilde{S}%
_{cc}^{\ast^{-1}}\left(  f_{c}-\left(  S_{rc}R_{c}\right)  ^{T}S_{rr}%
^{-1}f_{r}\right) \\
&  =M_{P-FETI-DP}r,
\end{align*}
where
\begin{align}
M_{P-FETI-DP}  &  =E_{r}S_{rr}^{-1}E_{r}^{T}+\label{eq:P-FETI-DP}\\
&  +\left(  E_{c}-E_{r}S_{rr}^{-1}S_{rc}R_{c}\right)  \widetilde{S}_{cc}%
^{\ast^{-1}}\left(  E_{c}^{T}-R_{c}^{T}S_{rc}^{T}S_{rr}^{-1}E_{r}^{T}\right)
\nonumber
\end{align}
is the associated preconditioner P-FETI-DP, same as~\cite[eq. (90)]%
{Fragakis-2003-MHP}.

\section{BDD Family of Methods}

\label{sec:BDD-methods}

We recall two primal preconditioners from the Balancing Domain Decomposition
(BDD) family by Mandel in~\cite{Mandel-1993-BDD}; namely the original BDD and
Balancing Domain Decomposition by Constraints (BDDC)\ introduced by
Dohrmann~\cite{Dohrmann-2003-PSC}.

\subsection{BDD}

The BDD is a Neumann-Neumann algorithm, cf., e.g.,~\cite{Dryja-1995-SMN},
with a simple coarse grid correction, introduced by Mandel~\cite%
{Mandel-1993-BDD}. The name of the preconditioner comes from an idea to
\emph{balance} the residual. We say that $v\in \widehat{W}$ is balanced if
\begin{equation*}
Z^{T}E^{T}v=0.
\end{equation*}%
Let us denote the \textquotedblleft balancing\textquotedblright\ operator as%
\begin{equation}
C=EZ,  \label{eq:C}
\end{equation}%
so the columns of $C$ are equal to the weighted sum of traces of the
subdomain zero energy modes. Next, let us denote by $S_{C}\widehat{S}$ the $%
\widehat{S}-orthogonal$ projection onto the range of $C$, so that
\begin{equation*}
S_{C}=C\left( C^{T}\widehat{S}C\right) ^{-1}C^{T},
\end{equation*}%
and by $P_{C}$ the complementary projection to $S_{C}\widehat{S}$,\ defined
as%
\begin{equation}
P_{C}=I-S_{C}\widehat{S}.  \label{eq:P_C}
\end{equation}%
The BDD preconditioner~\cite[Lemma~3.1]{Mandel-1993-BDD}, can be written in
our settings as%
\begin{eqnarray}
M_{BDD} &=&\left[ \left( I-S_{C}\widehat{S}\right) ES^{+}E^{T}\widehat{S}%
(I-S_{C}\widehat{S})+S_{C}\widehat{S}\right] \widehat{S}^{-1}  \notag \\
&=&\left[ \left( I-S_{C}\widehat{S}\right) ES^{+}E^{T}(\widehat{S}\widehat{S}%
^{-1}-\widehat{S}S_{C}\widehat{S}\widehat{S}^{-1})+S_{C}\widehat{S}\widehat{S%
}^{-1}\right]  \notag \\
&=&P_{C}ES^{+}E^{T}P_{C}^{T}+S_{C}  \label{eq:BDD}
\end{eqnarray}%
where $S_{C}$ serves as the coarse grid correction. See~\cite%
{Mandel-1993-BDD,Mandel-1996-BDD}, and~\cite{Fragakis-2003-MHP} for details.

\subsection{BDDC}

Following a similar path as Li and Widlund~\cite{Li-2006-FBB}, we will assume
that each constraint can be represented by an explicit degree of freedom and
that we can decompose the space $\widetilde{W}$\ as in (\ref{eq:tildeW-decomp}%
). We note that the original BDDC in~\cite{Dohrmann-2003-PSC,Mandel-2003-CBD}
is mathematically equivalent, but algorithmically it treats the corner coarse degrees of
freedom and edge in the definition of $\widetilde{W}$ in different ways. The
BDDC is the method of preconditioned conjugate gradients for the assembled
system~(\ref{eq:system}) with the preconditioner $M_{BDDC}$ defined by,
cf.~\cite[eq. (27)]{Li-2006-FBB},
\[
M_{BDDC}=T_{sub}+T_{0},
\]
where $T_{sub}=E_{r}S_{rr}^{-1}E_{r}^{T}$\ is the subdomain correction
obtained by solving independent problems on subdomains, and $T_{0}%
=E\Psi\left(  \Psi^{T}S\Psi\right)  ^{-1}\Psi^{T}E^{T}$\ is the coarse grid
correction. Here $\Psi$ are the coarse basis functions defined by energy
minimization,%
\[
\operatorname*{tr}\Psi^{T}S\Psi\rightarrow\min\text{.}%
\]
Since we assume that each constraint corresponds to an explicit degree of
freedom, the coarse basis functions$~\Psi$ can be easily determined via the
analogy\ to the discrete harmonic functions, discussed, e.g., in~\cite[Section
4.4]{Toselli-2005-DDM}; $\Psi$ are equal to $1$\ in the coarse degrees of
freedom and have energy minimal extension with respect to the remaining
degrees of freedom $u_{r}$, so they are precisely given as
\[
\Psi=\left(
\begin{array}
[c]{c}%
R_{c}\\
-S_{rr}^{-1}S_{rc}R_{c}%
\end{array}
\right)  .
\]
Then, we can compute%
\begin{align*}
\Psi^{T}S\Psi &  =\left(
\begin{array}
[c]{cc}%
R_{c}^{T} & -R_{c}^{T}S_{rc}^{T}S_{rr}^{-1}%
\end{array}
\right)  \left(
\begin{array}
[c]{cc}%
S_{cc} & S_{rc}^{T}\\
S_{rc} & S_{rr}%
\end{array}
\right)  \left(
\begin{array}
[c]{c}%
R_{c}\\
-S_{rr}^{-1}S_{rc}R_{c}%
\end{array}
\right) \\
&  =R_{c}^{T}S_{cc}R_{c}-R_{c}^{T}S_{rc}^{T}S_{rr}^{-1}S_{rc}R_{c}\\
&  =\widetilde{S}_{cc}-R_{c}^{T}S_{rc}^{T}S_{rr}^{-1}S_{rc}R_{c}\\
&  =\widetilde{S}_{cc}^{\ast},
\end{align*}
followed by%
\begin{align*}
&  E\Psi\left[  \Psi^{T}S\Psi\right]  ^{-1}\Psi^{T}E^{T}\\
&  =E\left(
\begin{array}
[c]{c}%
R_{c}\\
-S_{rr}^{-1}S_{rc}R_{c}%
\end{array}
\right)  \widetilde{S}_{cc}^{\ast^{-1}}\left(
\begin{array}
[c]{cc}%
R_{c}^{T} & -R_{c}^{T}S_{rc}^{T}S_{rr}^{-1}%
\end{array}
\right)  E^{T}\\
&  =\left(  E_{c}-E_{r}S_{rr}^{-1}S_{rc}R_{c}\right)  \widetilde{S}_{cc}%
^{\ast^{-1}}\left(  E_{c}^{T}-R_{c}^{T}S_{rc}^{T}S_{rr}^{-1}E_{r}^{T}\right)
.
\end{align*}
So, the BDDC\ preconditioner takes the form%
\begin{align}
M_{BDDC}  &  =E_{r}S_{rr}^{-1}E_{r}^{T}+\label{eq:BDDC}\\
&  +\left(  E_{c}-E_{r}S_{rr}^{-1}S_{rc}R_{c}\right)  \widetilde{S}_{cc}%
^{\ast^{-1}}\left(  E_{c}^{T}-R_{c}^{T}S_{rc}^{T}S_{rr}^{-1}E_{r}^{T}\right)
.\nonumber
\end{align}

\section{Connections of the Preconditioners}

\label{sec:connections}We review from~\cite[Section 8]{Fragakis-2003-MHP} that
a certain version of P-FETI-1 gives exactly the same algorithm as BDD. Next,
we state the equivalence of P-FETI-DP\ and BDDC preconditioners. Finally, we
translate the abstract proof relating the spectra of primal and dual
preconditioners~\cite[Theorem 4]{Fragakis-2007-FDD} in the case of FETI-1 and BDD.

\begin{theorem}
[{\cite[Section 8]{Fragakis-2003-MHP}}]\label{thm:FETI=BDD} If Q is chosen to
be the Dirichlet preconditioner, the P-FETI-1 and the BDD preconditioners\ are
the same.
\end{theorem}

\emph{Proof.} We will show that the P-FETI-1 in~(\ref{eq:P-FETI-1}) with
$Q=B_{D}SB_{D}^{T}$ is the same as the BDD in~(\ref{eq:BDD}). So, similarly as
in~\cite[pp. 3819-3820]{Fragakis-2003-MHP}, from (\ref{eq:H}) we get%
\begin{align*}
H  &  =I-B^{T}QG\left(  G^{T}QG\right)  ^{-1}Z^{T}\\
&  =I-B^{T}B_{D}SB_{D}^{T}BZ\left(  Z^{T}B^{T}B_{D}SB_{D}^{T}BZ\right)
^{-1}Z^{T}\\
&  =I-A_{R}\left(  Z^{T}A_{R}\right)  ^{-1}Z^{T},
\end{align*}
where
\[
A_{R}=B^{T}B_{D}SB_{D}^{T}BZ.
\]
Using (\ref{eq:jump+average}), definitions of $C$ in (\ref{eq:C}),
$\widehat{S}$\ in (\ref{eq:system}), and because $SZ=0$ by
(\ref{eq:rangeZ=nullS}),%
\begin{align*}
A_{R}  &  =\left(  I-E^{T}R^{T}\right)  S\left(  I-RE\right)  Z\\
&  =SZ-SREZ-E^{T}R^{T}SZ+E^{T}R^{T}SREZ\\
&  =SZ-SRC-E^{T}R^{T}SZ+E^{T}\widehat{S}C\\
&  =\left(  E^{T}\widehat{S}-SR\right)  C,
\end{align*}
and similarly%
\begin{align*}
Z^{T}A_{R}  &  =Z^{T}\left(  E^{T}\widehat{S}-SR\right)  C\\
&  =C^{T}\widehat{S}C-Z^{T}SREZ\\
&  =C^{T}\widehat{S}C.
\end{align*}
Using the two previous results, (\ref{eq:P_C})\ and symmetries of $\widehat
{S}$ and $S_{c}$, we get%
\begin{align*}
HE^{T}  &  =\left(  I-A_{R}\left(  Z^{T}A_{R}\right)  ^{-1}Z^{T}\right)
E^{T}\\
&  =E^{T}-A_{R}\left(  Z^{T}A_{R}\right)  ^{-1}Z^{T}E^{T}\\
&  =E^{T}-\left(  E^{T}\widehat{S}-SR\right)  C\left(  C^{T}\widehat
{S}C\right)  ^{-1}C^{T}\\
&  =E^{T}-\left(  E^{T}\widehat{S}-SR\right)  S_{C}\\
&  =E^{T}-E^{T}\widehat{S}S_{C}+SRS_{C}\\
&  =E^{T}\left(  I-\widehat{S}S_{C}\right)  +SRS_{C}\\
&  =E^{T}P_{C}^{T}+SRS_{C}.
\end{align*}
Next, the matrix $S_{C}$ satisfies the relation%
\begin{align*}
S_{C}R^{T}SS^{+}SRS_{C}  &  =S_{C}R^{T}SRS_{C}=S_{C}\widehat{S}S_{C}\\
&  =C\left(  C^{T}\widehat{S}C\right)  ^{-1}C^{T}\widehat{S}C\left(
C^{T}\widehat{S}C\right)  ^{-1}C^{T}\\
&  =C\left(  C^{T}\widehat{S}C\right)  ^{-1}C^{T}=S_{C}.
\end{align*}
Because by definition $P_{C}C=0$, using (\ref{eq:decomp-unity}) we get for
some $Y$ that
\begin{align*}
P_{C}ES^{+}SRS_{C}  &  =P_{C}E\left(  I+ZY\right)  RS_{C}\\
&  =P_{C}ERS_{C}+P_{C}EZYRS_{C}\\
&  =P_{C}S_{C}+P_{C}CYRS_{C}\\
&  =P_{C}S_{C}\\
&  =\left(  I-S_{C}\widehat{S}\right)  S_{C}\\
&  =S_{C}-S_{C}=0,
\end{align*}
and the same is true for the transpose, so $S_{C}R^{T}SS^{+}E^{T}P_{C}^{T}=0$.

Using these results, the P-FETI-1 preconditioner from~(\ref{eq:P-FETI-1}%
)\ becomes
\begin{align}
M_{P-FETI}  &  =EH^{T}S^{+}HE^{T}\nonumber\\
&  =\left(  S_{C}R^{T}S+P_{C}E\right)  S^{+}\left(  E^{T}P_{C}^{T}%
+SRS_{C}\right) \nonumber\\
&  =S_{C}R^{T}SS^{+}E^{T}P_{C}^{T}+S_{C}R^{T}SS^{+}SRS_{C}\nonumber\\
&  \quad+P_{C}ES^{+}E^{T}P_{C}^{T}+P_{C}ES^{+}SRS_{C}\nonumber\\
&  =P_{C}ES^{+}E^{T}P_{C}^{T}+S_{C}, \label{eq:P-FETI-1-Dirichlet}%
\end{align}
and we see that~(\ref{eq:P-FETI-1-Dirichlet}) is the same as the definition of
BDD in~(\ref{eq:BDD}). \qquad\endproof

\begin{theorem}
The P-FETI-DP and the BDDC preconditioners are the same.
\end{theorem}

\emph{Proof.} The claim follows directly comparing the definitions of both
preconditioners, P-FETI-DP in eq.~(\ref{eq:P-FETI-DP}) and the BDDC in
eq.~(\ref{eq:BDDC}). \qquad\endproof

\begin{corollary}
Comparing the preconditioner proposed by Cros~\cite[eq. 4.8]{Cros-2003-PSC}
with the definitions (\ref{eq:P-FETI-DP}) and (\ref{eq:BDDC}), it follows that
this preconditioner can be interpreted as either, P-FETI-DP\ or BDDC.
\end{corollary}

In the remaining, we will show the equality of eigenvalues of BDD and FETI-1,
with $Q$ being the Dirichlet preconditioner.

\begin{lemma}
\label{preconditioned-operators} The two preconditioned operators can be
written as%
\begin{align*}
M_{FETI}\mathcal{F}  &  =\left(  B_{D}SB_{D}^{T}\right)  \left(
B\widetilde{S}^{+}B^{T}\right)  ,\\
M_{BDD}\widehat{S}  &  =\left(  E\widetilde{S}^{+}E^{T}\right)  \left(
R^{T}SR\right)  ,
\end{align*}
where%
\[
\widetilde{S}^{+}=H^{T}S^{+}H.
\]

\end{lemma}

\emph{Proof.} First, $M_{FETI}=B_{D}SB_{D}^{T}$, which is the Dirichlet
preconditioner. From (\ref{eq:feti}), using the definition of $H$ by
(\ref{eq:H}), we get%
\begin{align*}
\mathcal{F}  &  =P^{T}FP\\
&  =P^{T}BS^{+}B^{T}P\\
&  =\left(  I-G(G^{T}QG)^{-1}G^{T}Q^{T}\right)  BS^{+}B^{T}\left(
I-QG(G^{T}QG)^{-1}G^{T}\right) \\
&  =\left(  B-BZ(G^{T}QG)^{-1}G^{T}Q^{T}B\right)  S^{+}\left(  B^{T}%
-B^{T}QG(G^{T}QG)^{-1}Z^{T}B^{T}\right) \\
&  =B\left(  I-Z(G^{T}QG)^{-1}G^{T}QB\right)  S^{+}\left(  I-B^{T}%
QG(G^{T}QG)^{-1}Z^{T}\right)  B^{T}\\
&  =BH^{T}S^{+}HB^{T}=B\widetilde{S}^{+}B^{T}\text{.}%
\end{align*}

Next, $\widehat{S}$ is defined by (\ref{eq:system}). By Theorem
\ref{thm:FETI=BDD}, we can\ use (\ref{eq:P-FETI-1}) for $M_{BDD}$\ to get%
\[
M_{BDD}=EH^{T}S^{+}HE^{T}=E\widetilde{S}^{+}E^{T}.\text{ }\qquad\endproof
\]

Before proceeding to the main result, we need to prove two technical Lemmas
relating the operators $S$ and $\widetilde{S}^{+}$. The first Lemma
establishes~\cite[Assumptions (13) and (22)]{Fragakis-2007-FDD} as well
as~\cite[Lemma~3]{Fragakis-2007-FDD} for FETI-1 and BDD.

\begin{lemma}
The operators $S$, $\widetilde{S}^{+}$ defined by (\ref{eq:block-operators})
and Theorem~\ref{preconditioned-operators}, resp., satisfy%
\begin{align}
\widetilde{S}^{+}SR  &  =R,\label{eq:id-1a}\\
\widetilde{S}^{+}S\widetilde{S}^{+}  &  =\widetilde{S}^{+}. \label{eq:id-1b}%
\end{align}
Moreover, the following relations are valid%
\begin{align}
B\widetilde{S}^{+}SR  &  =0,\label{eq:id-2}\\
\widetilde{S}^{+}B^{T}B_{D}S\widetilde{S}^{+}E^{T}  &  =0. \label{eq:id-3}%
\end{align}

\end{lemma}

\emph{Proof.} First, from (\ref{eq:rangeZ=nullS}) and symmetry of $S$ it
follows that%
\begin{align*}
HS  &  =\left(  I-B^{T}QG\left(  G^{T}QG\right)  ^{-1}Z^{T}\right)  S\\
&  =S-B^{T}QG\left(  G^{T}QG\right)  ^{-1}Z^{T}S=S.
\end{align*}
Using $H^{T}=I-Z\left(  G^{T}QG\right)  ^{-1}G^{T}QB$ we get%
\begin{align*}
H^{T}S^{+}S  &  =H^{T}\left(  I+ZY\right)  =H^{T}+H^{T}ZY\\
&  =H^{T}+\left[  I-Z\left(  G^{T}QG\right)  ^{-1}G^{T}QB\right]  ZY\\
&  =H^{T}+ZY-Z\left(  G^{T}QG\right)  ^{-1}G^{T}QGY\\
&  =H^{T}+ZY-ZY=H^{T},
\end{align*}
so
\[
\widetilde{S}^{+}S=H^{T}S^{+}HS=H^{T}S^{+}S=H^{T}.
\]
Finally, from previous and (\ref{eq:BR}), we get (\ref{eq:id-1a})\ as%
\[
\widetilde{S}^{+}SR=H^{T}R=\left(  I-Z\left(  G^{T}QG\right)  ^{-1}%
G^{T}QB\right)  R=R,
\]
and since $H^{T}$ is a projection, we immediately get also (\ref{eq:id-1b}) as%
\[
\widetilde{S}^{+}S\widetilde{S}^{+}=H^{T}\widetilde{S}^{+}=H^{T}H^{T}%
S^{+}H=\widetilde{S}^{+}.
\]
Next, (\ref{eq:id-2}) follows\ directly from (\ref{eq:id-1a}) noting
(\ref{eq:BR}).

Using (\ref{eq:id-1a})-(\ref{eq:id-1b}) and (\ref{eq:jump+average}%
)-(\ref{eq:B_DE^T}), we get (\ref{eq:id-3}) as
\begin{align*}
\widetilde{S}^{+}B^{T}B_{D}S\widetilde{S}^{+}E^{T}  &  =\widetilde{S}%
^{+}\left(  I-E^{T}R^{T}\right)  S\widetilde{S}^{+}E^{T}\\
&  =\widetilde{S}^{+}S\widetilde{S}^{+}E^{T}-\widetilde{S}^{+}E^{T}%
R^{T}S\widetilde{S}^{+}E^{T}\\
&  =\widetilde{S}^{+}E^{T}-\widetilde{S}^{+}E^{T}R^{T}E^{T}\\
&  =\widetilde{S}^{+}\left(  I-E^{T}R^{T}\right)  E^{T}\\
&  =\widetilde{S}^{+}B^{T}B_{D}E^{T}=0.\text{ }\qquad\endproof
\end{align*}

Next Lemma is a particular version of~\cite[Theorem 4]{Fragakis-2007-FDD} for
FETI-1 and BDD.

\begin{lemma}
\label{lem-precond-op-identities}The following identities are valid:
\begin{align*}
T_{D}\left(  M_{FETI}\mathcal{F}\right)   &  =\left(  M_{BDD}\widehat
{S}\right)  T_{D}, & T_{D}  &  =E\widetilde{S}^{+}B^{T},\\
T_{P}\left(  M_{BDD}\widehat{S}\right)   &  =\left(  M_{FETI}\mathcal{F}%
\right)  T_{P}, & T_{P}  &  =\left(  M_{FETI}\mathcal{F}\right)  B_{D}SR.
\end{align*}

\end{lemma}

\emph{Proof.} Using the transpose of (\ref{eq:id-3}) and (\ref{eq:id-2}), we
derive the first identity as%
\begin{align*}
T_{D}\left(  M_{FETI}\mathcal{F}\right)   &  =E\widetilde{S}^{+}B^{T}%
B_{D}SB_{D}^{T}B\widetilde{S}^{+}B^{T}\\
&  =E\widetilde{S}^{+}\left(  I-E^{T}R^{T}\right)  S\left(  I-RE\right)
\widetilde{S}^{+}B^{T}\\
&  =E\widetilde{S}^{+}S\left(  I-RE\right)  \widetilde{S}^{+}%
B^{T}-E\widetilde{S}^{+}E^{T}R^{T}S\widetilde{S}^{+}B^{T}\\
&  \quad+E\widetilde{S}^{+}E^{T}R^{T}SRE\widetilde{S}^{+}B^{T}\\
&  =E\widetilde{S}^{+}SB_{D}^{T}B\widetilde{S}^{+}B^{T}-E\widetilde{S}%
^{+}E^{T}R^{T}S\widetilde{S}^{+}B^{T}+\\
&  \quad+\left(  E\widetilde{S}^{+}E^{T}\right)  \left(  R^{T}SR\right)
T_{D}\\
&  =\left(  M_{BDD}\widehat{S}\right)  T_{D}\text{.}%
\end{align*}

Similarly, using (\ref{eq:id-3}) and (\ref{eq:id-2}), we derive the second
identity as%
\begin{align*}
T_{P}\left(  M_{BDD}\widehat{S}\right)   &  =\left(  M_{FETI}\mathcal{F}%
\right)  B_{D}SRE\widetilde{S}^{+}E^{T}R^{T}SR\\
&  =\left(  M_{FETI}\mathcal{F}\right)  B_{D}S\left(  I-B_{D}^{T}B\right)
\widetilde{S}^{+}\left(  I-B^{T}B_{D}\right)  SR\\
&  =\left(  M_{FETI}\mathcal{F}\right)  B_{D}S\widetilde{S}^{+}\left(
I-B^{T}B_{D}\right)  SR\\
&  \quad-\left(  M_{FETI}\mathcal{F}\right)  B_{D}SB_{D}^{T}B\widetilde{S}%
^{+}SR\\
&  \quad+\left(  M_{FETI}\mathcal{F}\right)  B_{D}SB_{D}^{T}B\widetilde{S}%
^{+}B^{T}B_{D}SR\\
&  =M_{FETI}B\widetilde{S}^{+}B^{T}B_{D}S\widetilde{S}^{+}E^{T}R^{T}SR\\
&  \quad-\left(  M_{FETI}\mathcal{F}\right)  B_{D}SB_{D}^{T}B\widetilde{S}%
^{+}SR\\
&  \quad+\left(  M_{FETI}\mathcal{F}\right)  \left(  B_{D}SB_{D}^{T}\right)
\left(  B\widetilde{S}^{+}B^{T}\right)  B_{D}SR\\
&  =\left(  M_{FETI}\mathcal{F}\right)  \left(  M_{FETI}\mathcal{F}\right)
B_{D}SR.\\
&  =\left(  M_{FETI}\mathcal{F}\right)  T_{P}\text{. }\qquad\endproof
\end{align*}

\begin{theorem}
Under the assumption of Lemma~\ref{lem-precond-op-identities}, the spectra of
the preconditioned operators $M_{BDD}\widehat{S}$\ and $M_{FETI-1}\mathcal{F}$
satisfy the relation%
\[
\sigma\left(  M_{BDD}\widehat{S}\right)  \setminus\left\{  1\right\}
=\sigma\left(  M_{FETI-1}\mathcal{F}\right)  \setminus\left\{  0,1\right\}  .
\]
Moreover, the multiplicity of any common eigenvalue $\lambda\neq0,1$ is
identical for the two preconditioned operators.
\end{theorem}

\emph{Proof.} Let $u_{D}$ be a (nonzero) eigenvector of the preconditioned
FETI-1 operator corresponding to the eigenvalue $\lambda_{D}$. Then, by
Lemma~\ref{lem-precond-op-identities}, we have%
\[
T_{D}\left(  M_{FETI-1}\mathcal{F}\right)  u_{D}=\left(  M_{BDD}\widehat
{S}\right)  T_{D}u_{D},
\]
so $T_{D}u_{D}$ is an eigenvector of the preconditioned BDD operator
corresponding to the eigenvalue $\lambda_{D}$, provided that $T_{D}u_{D}\neq
0$. So assume that $T_{D}u_{D}=0$. But then it is also true that
\begin{align*}
0  &  =B_{D}SR\left(  T_{D}u_{D}\right)  =B_{D}SRE\widetilde{S}^{+}B^{T}%
u_{D}\\
&  =B_{D}S\left(  I-B_{D}^{T}B\right)  \widetilde{S}^{+}B^{T}u_{D}%
=B_{D}S\widetilde{S}^{+}B^{T}u_{D}-B_{D}SB_{D}^{T}B\widetilde{S}^{+}B^{T}%
u_{D}\\
&  =B_{D}S\widetilde{S}^{+}B^{T}u_{D}-\left(  M_{FETI}\mathcal{F}\right)
u_{D}=B_{D}S\widetilde{S}^{+}B^{T}u_{D}-\lambda_{D}u_{D},
\end{align*}
so%
\[
B_{D}S\widetilde{S}^{+}B^{T}u_{D}=\lambda_{D}u_{D}.
\]
Note that, by~(\ref{eq:id-1a}) and (\ref{eq:BR}), we get
\begin{align*}
\left(  B_{D}S\widetilde{S}^{+}B^{T}\right)  ^{2}  &  =B_{D}S\widetilde{S}%
^{+}B^{T}B_{D}S\widetilde{S}^{+}B^{T}\\
&  =B_{D}S\widetilde{S}^{+}\left(  I-E^{T}R^{T}\right)  S\widetilde{S}%
^{+}B^{T}\\
&  =B_{D}S\widetilde{S}^{+}S\widetilde{S}^{+}B^{T}-B_{D}S\widetilde{S}%
^{+}E^{T}R^{T}S\widetilde{S}^{+}B^{T}\\
&  =B_{D}S\widetilde{S}^{+}B^{T}-B_{D}S\widetilde{S}^{+}E^{T}R^{T}B^{T}\\
&  =B_{D}S\widetilde{S}^{+}B^{T},
\end{align*}
so $B_{D}S\widetilde{S}^{+}B^{T}$ is a projection and therefore $\lambda
_{D}=0,1$.

Next, Let $u_{P}$ be a (nonzero) eigenvector of the preconditioned BDD
operator corresponding to the eigenvalue $\lambda_{P}$. Then, by
Lemma~\ref{lem-precond-op-identities}, we have%
\[
T_{P}\left(  M_{BDD}\widehat{S}\right)  =\left(  M_{FETI}\mathcal{F}\right)
T_{P},
\]
so $T_{P}u_{P}$ is an eigenvector of the preconditioned FETI-1 operator
corresponding to the eigenvalue $\lambda_{P}$, provided that $T_{P}u_{P}\neq
0$. So assume that $T_{P}u_{P}=0$. But then also using (\ref{eq:id-1a}) and
(\ref{eq:decomp-unity}), we get%
\begin{align*}
0  &  =T_{D}\left(  T_{P}u_{P}\right)  =T_{D}\left(  M_{FETI}\mathcal{F}%
\right)  B_{D}SRu_{P}\\
&  =\left(  M_{BDD}\widehat{S}\right)  T_{D}B_{D}SRu_{P}=\left(
M_{BDD}\widehat{S}\right)  E\widetilde{S}^{+}B^{T}B_{D}SRu_{P}\\
&  =M_{BDD}\widehat{S}E\widetilde{S}^{+}\left(  I-E^{T}R^{T}\right)  SRu_{P}\\
&  =M_{BDD}\widehat{S}E\widetilde{S}^{+}SRu_{P}-M_{BDD}\widehat{S}%
E\widetilde{S}^{+}E^{T}R^{T}SRu_{P}\\
&  =M_{BDD}\widehat{S}u_{P}-M_{BDD}\widehat{S}E\widetilde{S}^{+}E^{T}%
R^{T}SRu_{P}\\
&  =M_{BDD}\widehat{S}u_{P}-\left(  M_{BDD}\widehat{S}\right)  ^{2}u_{P},
\end{align*}
which is the same as
\[
\lambda_{P}u_{P}-\lambda_{P}^{2}u_{P}=\lambda_{P}\left(  1-\lambda_{P}\right)
u_{P}=0,
\]
and therefore $\lambda_{P}=0,1$.

Finally, let $\lambda\neq0,1$\ be an eigenvalue of the operator $M_{BDD}%
\widehat{S}$\ with the multiplicity $m$. From the previous arguments, the
eigenspace corresponding to $\lambda$\ is mapped by the operator $T_{P}$\ into
an eigenspace of $M_{FETI-1}\mathcal{F}$ and since this mapping is one-to-one,
the multiplicity of $\lambda$ corresponding to $M_{FETI-1}\mathcal{F}$ is
$n\geq m$. By the same argument, we can prove the opposite inequality and the
conclusion follows. \qquad\endproof

\bibliographystyle{siam}
\bibliography{longpaper}

\end{document}